\documentclass[11pt,a4paper]{article}
\usepackage{amsmath,mathtools}
\usepackage{amssymb}
\newtheorem{Thm}{Theorem}[section]
\newtheorem{Lem}[Thm]{Lemma}
\newtheorem{Cor}[Thm]{Corollary}
\newtheorem{Prop}[Thm]{Proposition}

\newtheorem{Rem}[Thm]{Remark}

\newtheorem{Rems}[Thm]{Remarks}

\newtheorem{Def}[Thm]{Definition}

\numberwithin{equation}{section}

\newcommand{\Z}{\mathbb{Z}}

\newcommand{\R}{\mathbb{R}}
\newcommand{\C}{\mathbb{C}}

\newcommand{\T}{\mathbb{T}}

\newcommand{\HH}{\mathcal{H}}
\newcommand{\bpr}{ \textit{Proof}: }

\newcommand{\eprsk}{~$\blacksquare$\medskip}

\newcommand{\eps}{\varepsilon}

\newcommand{\ra}{\rightarrow}
\newcommand{\SL}{\operatorname{SL}}
\newcommand{\PSL}{\operatorname{PSL}}
\newcommand{\HAP}{\operatorname{HAP}}
\newcommand{\SA}{\operatorname{SA}}
\newcommand{\Tr}{\operatorname{Tr}}
\newcommand{\la}{\langle}
\newcommand{\rg}{\rangle}
\newcommand{\Rr}{\mathcal R}
\newcommand{\Bb}{\mathcal B}

\providecommand{\AMS}{$\mathcal{A}$\kern-.1667em%
\lower.25em\hbox{$\mathcal{M}$}\kern-.125em$\mathcal{S}$}
\begin{document}% Do not forget this command!
%\issueinfo{xx}{yy}{2024} 
% \issueinfo{vv}{n}{yyyy}, vv is the volume, n is the number, yyyy is the year
% leave these information to fixed later by the editorial office
%\commby{Stefaan Vaes}% Editor is the name of the editor who accepted the article
%\pagespan{00}{www}{zzz}
% \pagespan{bbb}{eee} where bbb is the beginning page, eee is the ending page
\date{April 9, 2024}% This is the date of submission of the article; it
% will be fixed by the editorial office
%\revision{XXXX, 2024}% This is date(s) of revision of the manuscript; it
% will be fixed by the editorial office
\title{The linear $\SL_2(\Z)$-action on $\T^n$: ergodic and von Neumann algebraic aspects}% The title of the article

\author{PAUL JOLISSAINT and ALAIN VALETTE}% For multiple authors please use {\protec \and} sequence

%\begin{subjclass}
%\end{keywords}
\maketitle

\footnote{Primary 28D15; Secondary 46L36, 37A20, 43A07. Keywords: Ergodic p.m.p. action, amenable equivalence relation, group measure space construction, rigid inclusion Cartan subalgebra, type II$_1$ factor.}

\begin{abstract} The \textit{linear action} of $\SL_2(\R)$ on $\R^n$ corresponding to its unique irreducible representation induces an action $\SL_2(\Z)\curvearrowright \T^n$ for every $n\geq 2$ that factors through $\PSL_2(\Z)$ for $n$ odd. Thus,
setting $G_n=\SL_2(\Z)$ (resp. $G_n=\PSL_2(\Z)$) for $n$ even (resp. $n$ odd), $G_n\curvearrowright \T^n$ is free and ergodic, every ergodic sub-equivalence relation of the orbital equivalence relation is either amenable or rigid, and the fundamental group of the II$_1$ factor $N_n\coloneqq L^\infty(\T^n)\rtimes G_n$ is trivial. For $n$ even, $L^\infty(\T^n)\rtimes H$ is a maximal Haagerup subalgebra of $N_n$ for every suitable maximal amenable subgroup $H$ of $\SL_2(\Z)$.
\end{abstract}
%Section 1
\section{INTRODUCTION}

It is a classical fact that, for $n\geq 1$, the Lie group $\SL_2(\R)$ admits a unique irreducible representation on $\R^n$. Indeed, let $P_{n-1}(\R)$ be the space of homogeneous polynomials of degree $n-1$ in two variables $X,Y$ with real coefficients: we have a standard representation $\rho_{n-1}$ of $\SL_2(\R)$ on $P_{n-1}(\R)$ defined by 
$$
\Big(\rho_{n-1}\begin{pmatrix}
a & b \\c & d\end{pmatrix}\Big)(P)(X,Y)=P(aX+cY,bX+dY)
$$
for $P\in P_{n-1}(\R)$ and 
$\begin{pmatrix}
a & b \\c & d
\end{pmatrix}
\in \SL_2(\R)$. 
Now we may identify $P_{n-1}(\R)$ with $\R^n$ by means of the canonical basis $X^{n-1},X^{n-2}Y,\ldots,XY^{n-2},Y^{n-1}$. It is clear that the subgroup $P_{n-1}(\Z)$ of polynomials with integer coefficients is invariant under $\rho_{n-1}(\SL_2(\Z))$, so we get an $\SL_2(\Z)$-action on the $n$-dimensional torus $\T^n=P_{n-1}(\R)/P_{n-1}(\Z)$. If $n$ is odd, the representation $\rho_{n-1}$ factors through $\PSL_2(\R)$, so the action of $\SL_2(\Z)$ on $\T^n$ factors through $\PSL_2(\Z)$.

Define $G_n$ as $\SL_2(\Z)$ for $n$ even and as $\PSL_2(\Z)$ for $n$ odd. In this paper, by {\it linear action} on $\T^n$ we mean the action of $G_n$ on $\T^n$. Observe that the linear action preserves the normalized Lebesgue measure on $\T^n$, so it is a p.m.p. (probability measure-preserving) action. For $n=2$, we recover the well-studied action of $\SL_2(\Z)$ on the 2-torus. Our first result generalizes a well-known property of this case:

%Thm 1.1
\begin{Thm}\label{ergodic+free} 
Fix $n\geq 2$. Then the following statements hold true: 
\begin{enumerate}
\item [(i)] The linear action on $\T^n$ is ergodic and essentially free.
\item [(ii)] If $h\in \SL_2(\Z)$ is hyperbolic and $n$ is even, then the subgroup $\langle h\rangle$ generated by 
$h$ acts ergodically on $\T^n$.
\end{enumerate}
\end{Thm}

We observe that the second part of Theorem \ref{ergodic+free} only holds for $n$ even. Indeed we show:

%Prop 1.2
\begin{Prop}\label{nonergodic} For $n$ odd, no amenable subgroup of $\PSL_2(\Z)$ acts ergodically on $\T^n$.
\end{Prop}

A countable, ergodic, measure preserving equivalence relation $\mathcal{R}$ on a standard probability space $(X,\mu)$ is said to be {\it rigid} if the inclusion of $L^\infty(X,\mu)$ into the von Neumann algebra $L(\mathcal{R})$ associated with $\mathcal{R}$ (see \cite{FelMoo}), is rigid in the sense of Popa \cite{PopaAnnals}. In Theorem 0.1 of \cite{Ioana}, Ioana discovered a remarkable property of the orbital equivalence relation $\mathcal{S}_2$ induced by the linear action of $\SL_2(\Z)$ on $\T^2$, namely: every ergodic sub-equivalence relation of $\mathcal{S}_2$ is either amenable or rigid. We will refer to this phenomenon as {\it Ioana's dichotomy}, and we observe that it is a vast generalization of a result by Burger \cite{Burger} that says that, for a subgroup $H$ of $\SL_2(\Z)$, either $H$ is amenable (and so $\Z^2\rtimes H$ is amenable), or the pair $(\Z^2\rtimes H, \Z^2)$ has the relative property (T).  More examples of this dichotomy, for homogeneous spaces of certain Lie groups, were provided by Ioana and Shalom (Theorem C in \cite{IoSha}). We contribute to the list by showing:

%Thm 1.3
\begin{Thm}\label{IoanaDicho} For $n\geq 2$, let $\mathcal{S}_n$ be the orbital equivalence relation on $\T^n$ induced by the linear action. Then Ioana's dichotomy holds for $\mathcal{S}_n$.
\end{Thm}

In the final part of the paper, we consider von Neumann algebraic aspects, more precisely we deal with the von Neumann algebra $N_n\coloneqq L^\infty(\T^n)\rtimes G_n$ obtained via the group measure space construction. It is a type II$_1$ factor as a consequence of Theorem \ref{ergodic+free}. Denote by $L(G)$ the group von Neumann algebra of a countable group $G$; using Fourier series we may identify $N_n$ with $L(\Z^n\rtimes G_n)$. Now $G_n$ has the Haagerup property and the pair $(\Z^n\rtimes G_n,\Z^n)$ has the relative property (T) (by Proposition 3.1 in \cite{ValMax}), so the factors $N_n$ belong to Popa's $\mathcal{HT}$-class \cite{PopaAnnals}; hence by Section 6 and Corollary 8.2 of \cite{PopaAnnals} they inherit $\ell^2$-Betti numbers $\beta_k^{HT}(N_n)$ that coincide with those of the acting group $G_n$, and we have:

%Thm 1.4
\begin{Thm}\label{ThmTrivialFundGr}
For every $n\geq 2$, the $\ell^2$-Betti numbers of the factor $N_n$ satisfy $\beta_k^{HT}(N_n)=0$  for every $k\not=1$ and 
\[
\beta_1^{HT}(N_n)=
\begin{cases}
1/12 & \textrm{for\ even\ } n\\
1/6 & \textrm{for\ odd\ } n.
\end{cases}
\]
In particular, if $n$ is odd and $m$ is even, then $N_n$ and $N_m$ are non-isomorphic, and, for every $n$, the fundamental group $\mathcal{F}(N_n)$ is trivial.
\end{Thm}

 In Theorem 3.1 of \cite{JiSk}, Jiang and Skalski proved that, if $H$ is a maximal amenable subgroup of $\SL_2(\Z)$ containing some hyperbolic matrix (hence by Proposition 2.6 in \cite{ValMax}, the subgroup $H$ is isomorphic either to $\Z\times C_2$ or to $\Z\rtimes C_4$), then $L^\infty(\T^2)\rtimes H$ is a maximal Haagerup von Neumann subalgebra of $N_2$, where the Haagerup property for finite von Neumann algebras has been studied first in \cite{Cho} and then systematically by the first author in \cite{Jol}. Building on Theorem \ref{IoanaDicho}, we extend that result as follows:

%Thm 1.5
\begin{Thm}\label{MaxHaagerup} Fix $n\geq 2$. Then the following statements hold true:
\begin{enumerate}
\item [(1)] For any hyperfinite subfactor $P$ of $N_n$ which contains $L^\infty(\T^n)\cong L(\Z^n)$ as a Cartan subalgebra, every maximal Haagerup von Neumann subalgebra $P\subset Q\subset N_n$ is a hyperfinite subfactor of $N_n$. Moreover:
\item [(2)] If $n$ is even, let $R=L(\Z^n\rtimes H)$ where $H$ is a maximal amenable subgroup of $\SL_2(\Z)$ which contains a hyperbolic element. Then $R$ is a maximal Haagerup von Neumann subalgebra of $N_n$.
\item [(3)] If $n$ is odd, then $P$ and $Q$ in (1) are not of the form $L(K)$ for any subgroup $K$ of $\Z^n\rtimes\PSL_2(\Z)$. 
\end{enumerate}
\end{Thm}

The existence of a hyperfinite subfactor $P$ as in the preceeding theorem is a consequence of Proposition 3.6 of \cite{PopaInv}. 

Moreover, in a finite von Neumann algebra, every von Neumann subalgebra with the Haagerup property is contained in a maximal one: it is stated in Corollary 1.12 of \cite{JiSk} as a consequence of Lemma 1.11 of the same article. However, the latter only treats increasing unions of countably many subalgebras, which is not sufficient to apply Zorn's lemma, even if their arguments would work for nets instead of sequences; as the proof of the above mentioned lemma does not contain much details, we think that it is worth providing a proof of the existence of maximal Haagerup subalgebras in Proposition \ref{ExistHMax} below.

By a celebrated result of Connes-Feldman-Weiss \cite{CFW}, the sub-equi-valence relation $\Rr$ of $\mathcal{S}_n$ such that $P=L(\Rr)$ is induced by a $\Z$-action. Part 3 of Theorem \ref{MaxHaagerup} above shows that $\mathcal{R}$ will not be produced by any subgroup of $G_n$.

Section 2 contains the proofs of Theorem \ref{ergodic+free} and Proposition \ref{nonergodic}. Section 3 is about Ioana's dichotomy: the main new ingredient is the fact that the $\SL_2(\R)$-action induced on the projective space $\mathbb{P}^{n-1}(\R)$ by the representation $\rho_{n-1}$, is amenable (see Proposition \ref{TopolAmen}): for $n=2$, i.e. the case treated in \cite{Ioana}, it is clear that the group $\SL_2(\R)$ acts amenably on $\mathbb{P}^1(\R)$, as $\mathbb{P}^1(\R)$ identifies with $\SL_2(\R)/B$, where $B$ is the Borel subgroup, and $B$ is amenable (see Proposition 4.3.2 in \cite{Zim}). The final section 4 contains the proofs of Theorems \ref{ThmTrivialFundGr} and \ref{MaxHaagerup}.

%\par\vspace{-7mm}

%Section 2
\section{ERGODICITY}

{\emph{Proof of Theorem \ref{ergodic+free}:}} 
\begin{enumerate}
\item We start with essential freeness. Observe that $\rho_{n-1}(\SL_2(\R))\subset \SL_n(\R)$ (otherwise, composing $\rho_{n-1}$ with the determinant we would get a non-trivial homomorphism from $\SL_2(\R)$ to an abelian group); hence we have $\rho_{n-1}(\SL_2(\Z))\subset \SL_n(\R)\cap N_n(\Z)=\SL_n(\Z)$, so that the linear action on $\T^n$ is the restriction to $\rho_{n-1}(\SL_2(\Z))$ of the standard action of $\SL_n(\Z)$ on $\T^n$. The latter is known to be essentially free, see e.g. Lemma 5.2.4 in \cite{Zim} and the comment following it.

Concerning ergodicity, view $\T^n$ as a compact abelian group on which $\SL_2(\Z)$ acts by automorphisms. View $\Z^n$ as the Pontryagin dual of $\T^n$. By Proposition 1.5 in \cite{BekMay}, the linear action is ergodic if and only if every non-zero $\SL_2(\Z)$-orbit on $\Z^n$ is infinite. So assume by contradiction that some non-zero $v\in \Z^n$ has a finite orbit. Then the stabilizer $\Gamma$ of $v$ has finite index in $\SL_2(\Z)$, in particular $\Gamma$ is non-amenable. By Proposition 3.1 in \cite{ValMax}, the restriction of $\rho_{n-1}$ to $\Gamma$ is irreducible, so it cannot have non-zero invariant vectors,  and we reached a contradiction.

\item By Example 1.4.(iii) in \cite{BekMay}, it is enough to check that none of the complex eigenvalues of $\rho_{n-1}(h)$ is a root of unit. But for some $\lambda\in\R$ with $|\lambda|>1$, the matrix $h$ is conjugate in $\SL_2(\R)$ to 
$h_\lambda=
\begin{pmatrix}
\lambda & 0 \\
0 & 1/\lambda
\end{pmatrix}$. 
Now $\rho_{n-1}(h_\lambda)$ 
is diagonal in the basis $X^{n-1},X^{n-2}Y,\ldots,XY^{n-2},Y^{n-1}$ of $P_{n-1}(\R)$, with eigenvalues $\lambda^{n-1}, \lambda^{n-3},\ldots,\lambda^{-(n-3)},\lambda^{-(n-1)}$. Since $n$ is even, none of those eigenvalues is a root of 1. \eprsk

\end{enumerate}

{\emph{Proof of Proposition \ref{nonergodic}:}} Fix some odd $n$. Let $H$ be an amenable subgroup of $\PSL_2(\Z)$, we must show that $H$ has a finite orbit in $\Z^n\setminus\{0\}$. Clearly we may assume that $H$ is maximal amenable, so we may appeal to the classification in Lemma 2.5 of \cite{ValMax}, giving three cases to consider:
\begin{itemize}
\item $H$ is cyclic of order $3$, and there is nothing to prove.
\item $H$ is infinite cyclic consisting of parabolic elements: then $H$ is conjugate to the subgroup generated by 
$u=
\begin{pmatrix}
1 & 1 \\0 & 1
\end{pmatrix}$; but $\rho_{n-1}(u)$ fixes the monomial $X^{n-1}$.
\item $H$ is either infinite cyclic or infinite dihedral, with torsion-free part being generated by some hyperbolic matrix $h$, conjugate in $\PSL_2(\R)$ to a diagonal matrix 
$h_\lambda=\begin{pmatrix}
\lambda & 0 \\0 & 1/\lambda
\end{pmatrix}$ 
with $\lambda>1$. As in the proof of the second part of Theorem \ref{ergodic+free}, the eigenvalues of $\rho_{n-1}(h_\lambda)$ (hence also of $\rho_{n-1}(h)$) are $\lambda^{n-1}, \lambda^{n-3},\ldots,\lambda^{-(n-3)},\lambda^{-(n-1)}$, but now $1$ appears in the list, as $n$ is odd. So $\rho_{n-1}(h)$ admits in $\R^n$ an eigenvector $v$ with eigenvalue $1$. As $\rho_{n-1}(h)$ has integer coefficients, by Gaussian elimination we may assume $v$ to have rational coordinates, and by chasing denominators we may assume $v\in\Z^n$. Then the orbit of $v$ under $\rho_{n-1}(H)$ has at most two elements.
\eprsk
\end{itemize}

%Section 3
\section{IOANA'S DICHOTOMY}

%Ssection 3.1
\subsection{EQUIVALENCE RELATIONS ON STANDARD SPACES}

Let $(X,\Bb,\mu)$ be a standard probability space. A \textit{(countable) standard equivalence relation} on $X$ is a Borel subset $\Rr\subset X\times X$ which is an equivalence relation whose equivalence classes $[x]_\Rr\coloneqq \{y\in X\colon x\sim y\}$ are all at most countable. We denote by $p^i:X\times X\rightarrow X$ the projection defined by $p^i(x_1,x_2)=x_i$ for $i=1,2$; as they are countable-to-one maps, they imply the existence of $\sigma$-finite measures $\nu_i$, $i=1,2$, on the Borel subsets of $\Rr$ which are defined by
\[
\nu_i(C)\coloneqq \int_X |(p^i)^{-1}(x)\cap C|d\mu(x)
\]
for every such Borel set $C$. We assume henceforth that $\mu$ is \textit{invariant} under $\Rr$, which means that $\nu\coloneqq\nu_1=\nu_2$. Thus, for every integrable function $f$ on $\Rr$, one has
\[
\int_\Rr f(x,y)d\nu(x,y)=\int_X\Big(\sum_{x\sim y} f(x,y)\Big)d\mu(x)=
\int_X\Big(\sum_{y\sim x}f(x,y)\Big)d\mu(y).
\] 
Finally, the \textit{full group} of $\Rr$, denoted by $[\Rr]$, consists in all automorphisms $\theta$ of $(X,\Bb,\mu)$ such that $\theta(x)\in [x]_\Rr$ for almost every $x\in X$, and $\Rr$ is \textit{ergodic} if any $\Rr$-invariant Borel subset of $X$ is either null or co-null.

Natural examples of countable equivalence relations are given by p.m.p. actions of countable groups: if $\Gamma\curvearrowright (X,\mu)$ is such an action, then the {\it orbital equivalence relation} $\Rr_\Gamma\coloneqq\{(x,\gamma x)\colon x\in X, \gamma\in\Gamma\}$ is a countable equivalence relation which is ergodic if and only if the action of $\Gamma$ is. It turns out that, by Theorem 3 of \cite{FelMooI}, every equivalence relation is of this form, but the acting group is far from being unique.

Recall that $\Rr$ is \textit{hyperfinite} if there exists a Borel automorphism $T$ of $(X,\mu)$ such that, for (almost) every $x\in X$, the equivalence class $[x]_\Rr$ coincides with the orbit $\{T^kx\colon k\in \Z\}$ of $x$. It turns out that, by the main result of \cite{CFW}, hyperfiniteness of $\Rr$ is equivalent to its amenability, whose definition we recall from \cite{CFW}:

%Def 3.1
\begin{Def}\label{amenrel} With $X,\Rr$ as above: say that $\Rr$ is an \textit{amenable equivalence relation} if there is an assignment $x\mapsto \nu_x$ from $X$ to the state space of 
$\ell^\infty([
x]_\Rr)$ such that:
\begin{enumerate}
\item [(1)] the family $(\nu_x)_{x\in X}$ is $\Rr$-invariant, i.e. $\gamma_*(\nu_x)=\nu_y$ for every $\gamma=(y,x)\in\Rr$;
\item [(2)] for every $f\in L^\infty(\Rr,\nu)$, the function $x\mapsto \nu_x(f)$ is measurable.
\end{enumerate}
\end{Def}

%Ssection 3.2
\subsection{AMENABILITY OF THE LINEAR ACTION OF $\SL_2(\R)$ ON $\mathbb{P}^{n-1}(\R)$}

The proof of Ioana's main theorem in \cite{Ioana} uses in a crucial way the topological amenability of the action of $\SL_2(\Z)$ on $\mathbb{P}^1(\R)$ (in the sense of \cite{Anan}). Analogously, the proof of Theorem \ref{IoanaDicho} will rest on the topological amenability of the action of $\SL_2(\Z)$ on $\mathbb{P}^{n-1}(\R)$. 
Let us first recall the relevant definition.

\bigskip
Let $\Gamma$ be a countable group and $X$ be a locally compact space equipped with a left action of $\Gamma$ by homeomorphisms. We denote by $\mathrm{Prob}(\Gamma)$ the set of probability measures on $\Gamma$, equipped with the weak$^*$-topology, and with the natural action of $\Gamma$ inherited by the left translation action of $\Gamma$ on itself.

%Def 3.2
\begin{Def}\label{DefTopolAmenAction}
(\cite{Anan}, Definition 2.1) Let $(\Gamma,X)$ be as above. We say that the action is \textit{topologically amenable} if there exists a net $(m_i)_{i\in I}$ of continuous maps $x\mapsto m_i^x$ from $X$ to $\mathrm{Prob}(\Gamma)$ such that 
\[
\lim_i \Vert sm_i^x-m_i^{sx}\Vert_1=0
\]
uniformly on compact subsets of $\Gamma\times X$.
\end{Def}

We need the following stability properties of the class of topologically amen-able actions.

%Lemma 3.3
\begin{Lem}\label{StabAmen} Let $(\Gamma,X)$ be a topologically amenable action on a locally compact space $X$. 
\begin{enumerate}
\item If $Z$ is a closed $\Gamma$-invariant subspace of $X$, then the action of $\Gamma$ on $Z$ is topologically amenable.
\item Let $Y$ be any locally compact space endowed with a $\Gamma$-action. Then the diagonal action of $\Gamma$ on $X\times Y$ is topologically amenable.
\item Let $F$ be a finite group acting on $X$, such that the $F$-action commutes with the $\Gamma$-action, then the $\Gamma$-action on the orbit space $X/F$ is topologically amenable.
\item For every $n\geq 2$, the diagonal action of $\Gamma$ on the {\it $n$-fold symmetric product} $\Sigma^n X\coloneqq X^n/Sym(n)$ is topologically amenable.
\end{enumerate}
\end{Lem}

\bpr Let $(m_i)_{i\in I}$ be a net of continuous maps $X\rightarrow \mathrm{Prob}(\Gamma)$ as in Definition \ref{DefTopolAmenAction}.
\begin{enumerate} 
\item Restrict the maps $m_i$ to $Z$.
\item For $i\in I$, define $n_i:X\times Y\rightarrow \mathrm{Prob}(\Gamma)$ by $n_i^{(x,y)}\coloneqq m_i^x$ (for $(x,y)\in X\times Y$). Using the $\Gamma$-equivariance of the projection $X\times Y\rightarrow X$, the net $(n_i)_{i\in I}$ realizes the topological amenability of the action $(\Gamma,X\times Y)$.
\item Say that $F$ acts on $X$ on the right. For $i\in I$ and $f\in F$, define $m_{i,f}:X\rightarrow \mathrm{Prob}(\Gamma)$ by $m_{i,f}^x\coloneqq m_i^{xf}$ (for $x\in X$), and set $n_i=:\frac{1}{|F|}\sum_{f\in F}m_{i,f}$. Then the net $(n_i)_{i\in I}$ factors through $X/F$ and realizes the topological amenability of the action $(\Gamma,X/F)$.
\item Follows immediately from the two previous points.
\end{enumerate} \eprsk

It is apparently classical that $\mathbb{P}^n(\C)$ is homeomorphic to $\Sigma^n\mathbb{P}^1(\C)$, see e.g. Example 1.2 in \cite{Mostovoy}. Since we need an $\SL_2(\C)$-equivariant version of this homeomorphism, we elaborate on this. Denote by $P_n(\C)$ the space of homogeneous polynomials of degree $n$ in two variables $X,Y$ with complex coefficients: complexifying $\rho_n$ we get a representation of $\SL_2(\C)$ on $P_n(\C)$, identified with $\C^{n+1}$ through the canonical basis of monomials. Passing to the complex projective space $\mathbb{P}^n(\C)$, we get the linear action of $\SL_2(\C)$ on $\mathbb{P}^n(\C)$.

%Lemma 3.4
\begin{Lem}\label{homeomorphism} For every $n\geq 1$, there exists an $\SL_2(\C)$-equivariant homeomorphism $F:\Sigma^n\mathbb{P}^1(\C)\rightarrow \mathbb{P}^n(\C)$.
\end{Lem}

\bpr Consider the continuous, $\SL_2(\C)$-equivariant map:
$$
\begin{array}{cl}
(\C^2\setminus\{(0,0)\})^n &\rightarrow P_n(\C)\\
((a_1,b_1),\ldots,(a_n,b_n))&\mapsto P(X,Y)=\prod_{i=1}^n(a_iX+b_iY).
\end{array}
$$
It clearly descends to a continuous, $\SL_2(\C)$-equivariant map $F:\Sigma^n\mathbb{P}^1(\C)\rightarrow\mathbb{P}^n(\C)$, and it remains to check that $F$ is bijective. Let $P(X,Y)$ be a non-zero polynomial in $P_n(\C)$, write $P(X,Y)=\sum_{k=0}^n c_kX^{n-k}Y^k$, and let $\ell$ be the smallest index $k$ such that $c_k\neq 0$, so that 
$$
P(X,Y)=\sum_{k=\ell}^n c_kX^{n-k}Y^k=Y^n\sum_{k=\ell}^n c_k\Big(\frac{X}{Y}\Big)^{n-k}.
$$ 
Set $Z=\frac{X}{Y}$ and $Q(Z)=\sum_{k=\ell}^n c_kZ^{n-k}$; view $Q$ as a polynomial of degree $n-\ell$ in the polynomial ring $\C[Z]$; factoring $Q$ uniquely as a product of $n-\ell$ factors of degree 1, we get $Q(Z)=c_\ell\prod_{i=1}^{n-\ell}(Z-z_i)$, hence $P(X,Y)=c_\ell Y^\ell\prod_{i=1}^{n-\ell}(X-z_iY)$. So injectivity of $F$ follows from the unique factorization property in $\C[Z]$, and surjectivity follows from $\C$ being algebraically closed.
\eprsk

%Prop 3.5
\begin{Prop}\label{TopolAmen} For $K=\R,\C$, let $\Gamma$ be a discrete subgroup of $\SL_2(K)$. Then the linear action of $\Gamma$ on $\mathbb{P}^n(K)$ is topologically amenable.
\end{Prop} 

\bpr We start with $K=\C$. By Lemma \ref{homeomorphism} and part (iv) of Lemma \ref{StabAmen}, it is enough to prove the result for $n=1$. But this is classical; as $\mathbb{P}^1(\C)=\SL_2(\C)/B$ (where $B$ is the subgroup of upper triangular matrices) and $B$ is amenable, the action of $\Gamma$ on $\SL_2(\C)/B$ is topologically amenable (see Example 2.7(5) in \cite{Anan}).

For $K=\R$, since $\mathbb{P}^n(\R)$ is a closed $\SL_2(\R)$-invariant subspace of $\mathbb{P}^n(\C)$, the result follows from part (i) of Lemma \ref{StabAmen}.
\eprsk

%Ssection 3.3
\subsection{RIGIDITY FOR INCLUSIONS AND EQUIVALENCE RELATIONS}

In the present article, all von Neumann algebras are assumed to have separable preduals, and all considered tracial states are assumed to be normal and faithful; they will be also called traces for short. Let then $N$ be a finite von Neumann algebra equipped with a trace $\tau$, and $1\in B\subset N$ be a von Neumann subalgebra of $N$ with the same unit. Recall that $N$ identifies with a dense subspace $\hat N$ of the Hilbert space $L^2(N,\tau)$ so that $\langle \hat x,\hat y\rangle=\tau(y^*x)$ for all $x,y\in N$; recall also that the map $J:\hat x\mapsto \widehat x^*$ extends to an antilinear isometry on $L^2(N,\tau)$ such that the operator $Jy^*J\in B(L^2(N,\tau))$ acts on the right: $Jy^*J\hat x=\widehat{xy}$ for all $x,y\in N$. 

We will consider below the following particular class of von Neumann subalgebras of $N$: a \textit{Cartan subalgebra} $A$ of $N$ is a maximal abelian von Neumann subalgebra of $N$ whose normalizer $\mathcal{N}_N(A)\coloneqq \{u\in U(N)\colon uAu^*=A\}$ generates $N$.

Important von Neumann algebras that will be considered here are the associated von Neumann algebras of groups: let $G$ be a countable group; its von Neumann algebra $L(G)$ is generated by the left regular representation $g\mapsto u_g$ on $\ell^2(G)$. It is finite, and it is equipped with the natural trace $\tau$ defined by $\tau(x)=\la x\delta_e,\delta_e\rg$ for all $x\in L(G)$. Every $x\in L(G)$ has a unique "Fourier series" expansion $\sum_g x(g)u_g$, where $x(g)=\tau(xu_{g^{-1}})$ and $\sum_g|x(g)|^2=\Vert x\Vert_2^2$. If $H<G$ is a subgroup, then $L(H)$ identifies in a natural way with a von Neumann subalgebra of $L(G)$, its elements $y$ being characterized by the fact that $y(g)=0$ for all $g\in G\setminus H$. Recall also that $L(G)$ is a factor if and only if $G$ is icc, i.e. it has infinite conjugacy classes, except the class of $e$.

More generally, assume that there is a trace-preserving action $\sigma$ of a group $G$ on a finite von Neumann algebra $(B,\tau)$. Then the associated crossed product is a finite von Neumann algebra acting on $L^2(B)\otimes \ell^2(G)$; it is generated by a copy of $B$ and a unitary representation $g\mapsto u_g$ of $G$ in such a way that $u$ is quasi-equivalent to the regular representation and that $u_gbu_g^*=\sigma_g(b)$ for all $g\in G$ and $b\in B$. It is denoted by $B\rtimes_\sigma G$.

Let us now recall the definition of von Neumann algebras associated to equivalence relations as in \cite{FelMoo}.

Let $\Rr$ be a standard equivalence relation and let $s$ be a normalized $2$-cocycle on $\Rr$: the latter is a Borel $\T$-valued map defined on 
\[
\Rr^2\coloneqq\{(x,y,z)\in X^3 \colon x\sim y\sim z\}
\] 
which has the following properties:
\begin{enumerate}
\item [(a)] $s(x,y,z)=1$ as soon as $|\{x,y,z\}|\leq 2$;
\item [(b)] for all $x_1\sim x_2\sim x_3\sim x_4$,
\[
s(x_1,x_3,x_4)s(x_1,x_2,x_3)=s(x_1,x_2,x_4)s(x_2,x_3,x_4).
\]
\end{enumerate}
The von Neumann algebra $L(\Rr,s)$ acts on $\HH=L^2(\Rr,\nu)$ and it is generated by the following two families of operators:
\begin{enumerate}
\item [(1)] We associate to $\phi\in [\Rr]$ the operator $u_\phi$ on $\HH$ defined by
\[
(u_\phi\xi)(x,z)=\xi(\phi^{-1}(x),z)s(x,\phi^{-1}(x),z)=\sum_{y\sim x}
\delta_{y,\phi^{-1}(x)}\xi(y,z)s(x,y,z) ;
\]
\item [(2)] next, we associate to every $a\in L^\infty(X,\mu)$ the multiplication operator on $\HH$, still denoted by $a$, 
which is defined by
\[
(a\xi)(x,z)=a(x)\xi(x,z)(=a(x)\xi(x,z)\underbrace{s(x,x,z)}_{=1}).
\]
\end{enumerate}
These two families of operators satisfy the following relations which are straightforward consequences of the definitions and whose proofs are left to the reader:

%Lemma 3.6
\begin{Lem}\label{ProprL}
Retain the above definitions. Then
\begin{enumerate}
\item [(A)] If $\theta,\phi\in [\Rr]$, then $u_\theta u_\phi=\sigma_{\theta,\phi}\cdot u_{\theta\phi}$ where $\sigma_{\theta,\phi}\in L^\infty(X,\mu)$ is the function
\[
\sigma_{\theta,\phi}(x)=s(x,\theta^{-1}(x),(\theta\phi)^{-1}(x)).
\]
In particular $u_{\theta^{-1}}u_\theta=1$ since $\sigma_{\theta^{-1},\theta}(x)=s(x,\theta(x),x)=1$ for every $x$.
\item [(B)] For every $\theta\in [\Rr]$, one has $u_\theta^*=u_{\theta^{-1}}$. Thus $u_\theta$ is a unitary operator.
\item [(C)] For all $a\in L^\infty(X,\mu)$ and $\phi\in [\Rr]$, one has
$u_\phi a u_\phi^*=a\circ \phi^{-1}$. \eprsk
\end{enumerate}
\end{Lem}
Thus, the set of linear combinations of the form $a u_\phi$ is a $*$-subalgebra of $B(\HH)$, and by definition, $L(\Rr,s)$ is its weak-operator closure, or its bicommutant, which is the same. It is a finite von Neumann algebra equipped with the natural trace
\[
\tau(y)=\langle y1_\Delta,1_\Delta\rangle \quad (y\in L(\Rr,s))
\]
where $\Delta=\{(x,x)\colon x\in X\}$ is the diagonal of $\Rr$. 

When $s=1$ is the trivial cocycle, we write $L(\Rr)$ instead of $L(\Rr,1)$.

One checks that $A=L^\infty(X,\mu)\subset N=L(\Rr,s)$ is a Cartan subalgebra, and that
$\Rr$ is ergodic if and only if $L(\Rr,s)$ is a II$_1$ factor. Conversely, Theorem 1 of \cite{FelMoo} states that if $A\subset N$ is a Cartan subalgebra of a finite von Neumann algebra with separable predual, then there exists an essentially unique pair $(\Rr,s)$ on a standard probability space $(X,\mu)$ such that $A\cong L^\infty(X)$ and $N=L(\Rr,s)$. (In fact, $N$ need not be finite, but we are only interested in this case here.) 

For future use, it is worth noting that, if $A\subset N$ is a Cartan subalgebra and if $A\subset Q\subset N$ is an intermediate subfactor, then $A$ is still a Cartan subalgebra of $Q$, by the comments preceeding Proposition 2.2 of \cite{PopaCartan}, and that there is a sub-equivalence relation $\mathcal{Q}$ of $\Rr$ such that $Q=L(\mathcal{Q},s|_\mathcal{Q})$, by Proposition 3.4 in \cite{AOI}.

Consider a p.m.p. ergodic and essentially free action $\Gamma\curvearrowright (X,\mu)$. Then, as explained in Subsection 3.1, $\Rr_\Gamma$ is an ergodic equivalence relation and the inclusion $L^\infty(X)\subset L(\Rr_\Gamma)$ is naturally isomorphic to the inclusion $L^\infty(X)\subset L^\infty(X)\rtimes\Gamma$. See for instance Section 1.2 in \cite{Ioana}.

%Remarks 3.7
\begin{Rems}\label{RemCrossedP}
(1) Observe that if $B$ is abelian and if the action is essentially free and ergodic, then it is a Cartan subalgebra of $B\rtimes_\sigma G$.\\
(2) If $\Gamma$ is a group acting on $H$ by automorphisms, then the action extends in a natural way to $L(H)$ and the algebra $L(H\rtimes \Gamma)$ is naturally isomorphic to the crossed product $L(H)\rtimes \Gamma$.
\end{Rems}

We are now ready to define the relative rigidity property of pairs $B\subset N$.
In order to do that, 
let us recall first that a Hilbert space $\HH$ is a \textit{Hilbert $N$-bimodule} if it admits commuting left and right normal actions. A vector $\xi\in\HH$ is \textit{tracial} if $\la x\xi,\xi\rg=\la \xi x,\xi\rg=\tau(x)$ for all $x\in N$, and it is \textit{$B$-central} if $b\xi=\xi b$ for all $b\in B$. A pair $(\HH,\xi)$ where $\xi\in \HH$ is a unit vector is called a \textit{pointed Hilbert $N$-bimodule}. 

The notion of rigidity of inclusions is due to Popa (\cite{PopaAnnals}, Sections 4 and 5); however, for practical reasons, we recall Definition 1.3.1 of \cite{Ioana}:

%Def 3.8
\begin{Def}\label{DefRigid}
The inclusion $B\subset N$ is \textit{rigid} (or has \textit{relative property (T)}) if, for every $\eps>0$, there exists a finite set $F\subset N$ and $\delta>0$ such that, if $(\HH,\xi)$ is a pointed Hilbert $N$-bimodule where $\xi$ is a tracial vector verifying 
\[
\max_{x\in F}\Vert x\xi-\xi x\Vert\leq \delta,
\]
then there exists a $B$-central vector $\eta\in\HH$ with $\Vert \eta-\xi\Vert\leq\eps$.
\end{Def}

%Remark 3.9
\begin{Rem}\label{RemRigid}
It is worth being mentioned that Definition \ref{DefRigid} is in fact independent of the chosen trace on $N$. Furthermore, if $B\subset N$ is a Cartan subalgebra of $N$, it follows from Corollary 1.3.3 of \cite{Ioana} that the required control on the distance between the $B$-central vector $\eta$ and the tracial vector $\xi$ can be omitted: the inclusion is rigid if and only if there exist $F\subset N$ finite and $\delta>0$ such that any pointed Hilbert $N$-bimodule $(\HH,\xi)$, where $\xi$ is tracial and $\Vert x\xi-\xi x\Vert\leq \delta$ for every $x\in F$, contains a non-zero $B$-central vector.
\end{Rem}

Important examples of rigid inclusions are provided by discrete groups:

%Prop 3.10
\begin{Prop}\label{ExRigid}
(Corollary 5.2, \cite{PopaAnnals}) Let $G$ be a countable group and let $H\subset G$ be a subgroup of $G$. Denote by $L(H)\subset L(G)$ their associated von Neumann algebras. Then the inclusion $L(H)\subset L(G)$ is rigid if and only if the pair $(G,H)$ has the relative property (T). \eprsk
\end{Prop}

Let us now recall Popa's definition of rigid equivalence relations:

%Def 3.11
\begin{Def}\label{DefRelRig}
(Definition 5.10.1, \cite{PopaAnnals}) Let $(X,\mu)$ be a standard probability
 space and let $\Rr$ be a countable, ergodic, measure preserving equivalence relation on $X$. Then $\Rr$ is \textit{rigid} if its associated inclusion $L^\infty(X)\subset L(\Rr)$ is rigid. Also, we say that an essentially free, ergodic, p.m.p. action $\Gamma\curvearrowright (X,\mu)$ is \textit{rigid} if the inclusion $L^\infty(X)\subset L^\infty(X)\rtimes\Gamma$ is rigid. Equivalently, the action $\Gamma\curvearrowright (X,\mu)$ is rigid if and only if the associated equivalence relation $\Rr_\Gamma$ is rigid.
\end{Def}

%Def 3.12
\begin{Def}\label{DefDicho}
Let $\Gamma\curvearrowright (X,\mu)$ be an ergodic and free p.m.p action. Then the action satisfies \textit{Ioana's dichotomy} if every ergodic sub-equivalence relation $\Rr\subset \Rr_\Gamma$ is either amenable (in the sense of Definition \ref{amenrel}) or rigid.
\end{Def}

 One of the main goals of the present article is to provide examples of Haage-rup maximal subalgebras that generalize Theorem 3.1 of \cite{JiSk}, as stated in Theorem \ref{MaxHaagerup}. Thus, let us consider the following more general setting: let $\Gamma\curvearrowright (X,\mu)$ be an ergodic and essentially free action so that $A=L^\infty(X)\subset N=A\rtimes\Gamma$ is the inclusion of the Cartan subalgebra $A$ in the type II$_1$ factor $N$. Assume furthermore that the action satisfies Ioana's dichotomy,  
and let $A\subset P$ be a subfactor which is maximal amenable in $N$. We intend to prove that $P$ is a maximal Haagerup subalgebra of $N$. In order to prove this, we consider an intermediate von Neumann algebra $P\subset Q\subset N$ which is assumed to have the Haagerup property, and we have to prove that $Q=P$. 
As a first step, as $A$ is maximal abelian in $N$, we observe that
\[
Q'\cap N\subset P'\cap N \subset A'\cap N=A\subset P,
\]
and as $P$ is a factor, the same holds for $Q$. Hence, by Theorem 1 of \cite{FelMoo}, there exists an ergodic equivalence relation $\Rr$ on $X$ and a Borel $2$-cocycle $s$ on the relation $\Rr$ such that $Q=L(\Rr,s)$, as defined after Lemma \ref{ProprL} \footnote{The authors of \cite{JiSk} appear to have overlooked the possibility that $s$ is non-trivial in the proof of their Theorem 3.1.}. 
In order to prove the equality $Q=P$, we need to prove that $\Rr$ is amenable, using Ioana's criterion of rigidity for equivalence relations as in Section 2 of \cite{Ioana}, but after having proved that the potential presence of the $2$-cocycle $s$ has no effect on the conclusion. In other words, we need to adapt Lemma 2.1 and Proposition 2.2 of \cite{Ioana}.

In order to do that, we denote by $B(X)$ the algebra of all complex-valued, bounded Borel functions on $X$. For two functions $f_1,f_2:X\ra\C$, the function $f_1\otimes f_2:X\times X\ra\C$ is defined by $(f_1\otimes f_2)(x_1,x_2)=f_1(x_1)f_2(x_2)$. Given two (regular) Borel probability measures $\mu$ and $\nu$ on $X$, the norm $\Vert\mu-\nu\Vert$ is defined by
\[
\Vert\mu-\nu\Vert=\sup\Big\{\Big|\int_Xfd\mu-\int_Xfd\nu\Big|\colon f\in B(X), \Vert f\Vert_\infty\leq 1\Big\}.
\]
Finally, observe that the quotient map $B(X)\ra L^\infty(X,\mu)$ makes any $L^\infty(X,\mu)$-bimodule a $B(X)$-bimodule.

%lem 3.13
\begin{Lem}\label{LemRigid}
Let $\Rr$ be a standard equivalence relation as above on $(X,\mu)$ so that $\mu$ is invariant under $\Rr$, and let $s$ be a Borel $2$-cocycle as above. Let $(\HH,\xi)$ be a pointed Hilbert $L(\Rr,s)$-bimodule. Then there exists a probability measure $\nu$ on $X\times X$ such that
\begin{enumerate}
\item [(i)] $\int_{X\times X}(f_1\otimes f_2)d\nu=\langle f_1\xi f_2,\xi\rangle$, $\forall f_1,f_2\in B(X)$.
\item [(ii)] $\Vert (\theta\otimes \theta)_*\nu-\nu\Vert\leq 2\Vert u_\theta \xi -\xi u_\theta\Vert$, $\forall \theta\in [\Rr]$.
\item [(iii)] If $\xi$ is tracial, then $p_*^i\nu=\mu$ for $i=1,2$.
\item [(iv)] If $\HH$ has no non-zero $L^\infty(X,\mu)$-central vector, then $\nu(\Delta)=0$.
\end{enumerate}
\end{Lem}
\bpr The proof is the same as that of Lemma 2.1 in \cite{Ioana} since the only place where $s$ occurs is in the definition of the unitary $u_\theta$ and, due to the latter, is used just in relation (C) of Lemma \ref{ProprL}, namely $u_\theta au_\theta^*=a\circ\theta^{-1}$ for $a\in L^\infty(X)$.
 \eprsk

%Prop 3.14
\begin{Prop}\label{PropRigid}
Let $\Rr$ be as in Lemma \ref{LemRigid}. Suppose moreover that $\Rr$ is ergodic. Assume that there is no sequence of probability measures $(\nu_k)_{k\geq 1}$ on $X\times X$ such that $\nu_k(\Delta)=0$ and $p_*^i\nu_k=\mu$ for all $i=1,2$ and $k$,
\begin{enumerate}
\item [(i)] $\lim_{k\to\infty}\int_{X\times X} (f_1\otimes f_2)d\nu_k=\int_X f_1f_2 d\mu$, $\forall f_1,f_2\in B(X)$, and
\item [(ii)] $\lim_{k\to\infty} \Vert (\theta\otimes \theta)_*\nu_k-\nu_k\Vert=0$, $\forall \theta\in [\Rr]$.
\end{enumerate}
Then, for every Borel $2$-cocycle $s$, the inclusion $L^\infty(X,\mu)\subset L(\Rr,s)$ is rigid.
\end{Prop}
\bpr Assume that there exists a Borel $2$-cocycle on $\Rr$ such that the inclusion $L^\infty(X,\mu)\subset L(\Rr,s)$ is not rigid. Using Definition \ref{DefRigid} and Remark \ref{RemRigid}, we can find a sequence $(\HH_k,\xi_k)$ of pointed Hilbert $L(\Rr,s)$-bimodules such that 
\[
\lim_{k\to\infty}\Vert x\xi_k-\xi_k x\Vert=0
\] 
for all $x\in L(\Rr,s)$, $\xi_k$ is tracial, and $\HH_k$ has no non-zero $L^\infty(X,\mu)$-central vector, for all $k\geq 1$. Let then $\nu_k$ be the probability measure associated to $(\HH_k,\xi_k)$ by Lemma \ref{LemRigid}. Then we conclude exactly as in the proof of Proposition 2.2 in \cite{Ioana}. 
 \eprsk

%Ssection 3.4
\subsection{PROOF OF THEOREM \ref{IoanaDicho}}

In analogy with Theorem 3.1 of \cite{Ioana}, we will prove:

%Thm 3.15
\begin{Thm}\label{IoanaMain} For $n\geq 2$, let $\mathcal{S}_n$ denote the orbital equivalence relation induced by the linear action on $\T^n$. Let $\mathcal{R}$ be an ergodic subequivalence relation of $\mathcal{S}_n$. If $\mathcal{R}$ is not amenable, then there does not exist a sequence of probability measures $(\nu_k)_{k\geq 1}$ on $\T^n\times\T^n$ such that $\nu_k(\Delta)=0$ for all $k$ and
\begin{enumerate}
\item[(i)] $\lim_{k\rightarrow\infty}\int_{\T^n\times\T^n}(f_1\otimes f_2)\,d\nu_k=\int_{\T^n}f_1f_2\,d\lambda_n$
for every $f_1,f_2\in B(\T^n)$; here $d\lambda_n$ denotes the normalized Lebesgue measure on $\T^n$;
\item[(ii)] $\lim_{k\rightarrow\infty}\Big|\int_{\T^n\times\T^n}(f\circ(\theta\times\theta))\,d\nu_k - \int_{\T^n\times\T^n} f\,d\nu_k\Big|=0$ for every $f\in B(\T^n\times\T^n)$ and $\theta\in[\mathcal{R}]$.
\end{enumerate}
\end{Thm}

The following consequence follows immediately by combining Theorem \ref{IoanaMain} with Proposition \ref{PropRigid}; it strengthens Theorem \ref{IoanaDicho}.

%Cor 3.16
\begin{Cor}\label{DichoWithCoc} For $n\geq 2$, let $\mathcal{R}$ be an ergodic sub-equivalence relation of $\mathcal{S}_n$. Then either $\mathcal{R}$ is amenable or, for every Borel $2$-cocycle $s$, the inclusion $L^\infty(X,\mu)\subset L(\Rr,s)$ is rigid. \eprsk
\end{Cor}

{\emph{Proof of Theorem \ref{IoanaMain}:}} Working by contrapositive, we assume that there exist probability measures $\nu_k$ on $\T^n\times \T^n$, with $\nu_k(\Delta)=0$ for every $k$, and satisfying (i) and (ii), and our aim is to prove that $\mathcal{R}$ is amenable. Now the proof is copying the 4 steps of the proof of Theorem 3.1 in \cite{Ioana}, replacing the $\SL_2(\Z)$-action on $\T^2$ by the $G_n$-action on $\T^n$, so we will merely indicate the non-obvious changes.

\vspace{3mm}

{\bf Step 1} goes through.

\vspace{3mm}

{\bf Notations:} We replace Ioana's maps $\sigma, r,\chi,q,p,\pi$ by the following ones:
\begin{itemize}
\item $\sigma:\R^n\setminus\{0\}\rightarrow \mathbb{P}^{n-1}(\R)$ is the canonical map, and $r=id\times\sigma:\T^n\times(\R^n\setminus\{0\})\rightarrow \T^n\times\mathbb{P}^{n-1}(\R)$.
\item $\chi:[-\frac{1}{2},\frac{1}{2}[^n\rightarrow\T^n$ is the continuous bijection given by $\chi(x_1,\ldots,x_n)=(e^{2\pi ix_1},\ldots,e^{2\pi ix_n})$; let $\rho:\T^n\setminus\{(1,\ldots,1)\}\rightarrow \R^n\setminus\{0\}$ denote the inverse of $\chi$ restricted to $\T^n\setminus\{(1,\ldots,1)\}$. Set $q=id\times\rho:\T^n\times(\T^n\setminus\{(1,\ldots,1)\})\rightarrow \T^{n}\times(\R^n\setminus\{0\})$.
\item $p:(\T^n\times\T^n)\setminus\Delta\rightarrow\T^n\times(\T^n\setminus\{(1,\ldots,1)\})$ is the homeomorphism given by $p(x,y)=(x,x^{-1}y)$.
\item Finally $$\pi=r\circ q\circ p:(\T^n\times\T^n)\setminus\Delta\rightarrow \T^n\times\mathbb{P}^{n-1}(\R):(x,y)\to(x,(\sigma\circ\rho)(x^{-1}y)).$$
\end{itemize}
Remark 3.3 of \cite{Ioana}, on $G_n$-equivariance of $\sigma, r$ and $p$, and limited equivariance of $\rho$, then goes through without change. 

\vspace{3mm}

{\bf Steps 2 and 3} go through without further change.

\vspace{3mm}

{\bf Step 4:} We must replace the topological amenability of the $\SL_2(\Z)$-action on $\mathbb{P}^1(\R)$ by the topological amenability of the $\SL_2(\Z)$-action on $\mathbb{P}^{n-1}(\R)$, which is Proposition \ref{TopolAmen}. Once this is done, the rest of the proof of Step 4 in \cite{Ioana} proceeds without change.
\eprsk

%Section 4
\section{VON NEUMANN ALGEBRAIC APPLICATIONS}

%Ssection 4.1
\subsection{RELATIVE PROPERTY H AND HT CARTAN SUBALGEBRAS}

As in Section 3, let $N$ be a finite von Neumann algebra with separable predual equipped with a trace $\tau$.  
For a von Neumann subalgebra $B\subset N$, let $E_B$ be the unique $\tau$-preserving conditional expectation from $N$ onto $B$. It extends to the orthogonal projection from $L^2(N,\tau)$ onto $L^2(B,\tau)$. More generally, if $\phi:N\ra N$ is a completely positive map such that $\tau\circ\phi\leq \tau$, then it is normal and it extends to a bounded operator $T_\phi$ on $L^2(N,\tau)$ characterized by $T_\phi(\hat x)=\widehat{\phi(x)}$ for every $x\in N$. Notice that $\phi$ is $B$-bimodular (i.e. $\phi(b_1xb_2)=b_1\phi(x)b_2$ for all $b_1,b_2\in B$ and $x\in N$) if and only if $T_\phi\in B'\cap JB'J$. The von Neumann algebra $JB'J$ coincides with the so-called \textit{basic construction} $\langle N,e_B\rangle$ which is generated by $N$ and the projection $e_B$ (see for instance Subsection 1.3.1 of \cite{PopaAnnals}). As its commutant $JBJ$ is a finite von Neumann algebra, $\langle N,e_B\rangle$ is semifinite; we then denote by $\mathcal{J}(\langle N,e_B\rangle)$ the norm-closed, two-sided ideal in 
$\langle N,e_B\rangle$ generated by its finite projections. If $B=\C1$, then $\la N,e_B\rg=JB'J$ coincides with the algebra of all bounded operators on $L^2(N,\tau)$ and the ideal $\mathcal{J}(\langle N,e_B\rangle)$ is the ideal of all compact operators $\mathcal{K}(L^2(N,\tau))$.

The following definition is due to Popa \cite{PopaAnnals}; notice that the Haagerup property considered in \cite{Jol} corresponds to the case $B=\C1$.

%Def 4.1
\begin{Def}\label{DefRelH}
Let $B\subset N$ be von Neumann algebras as above. Then $N$ has \textit{property H relative to} $B$ if there exists a sequence $(\phi_n)$ of normal, completely positive $B$-bimodular maps on $N$ satisfying the following conditions:
\begin{enumerate}
\item [(1)] $\phi_n(1)\leq 1$ and $\tau\circ \phi_n\leq \tau$ for every $n$;
\item [(2)] $(T_{\phi_n})\subset\mathcal{J}(\langle N,e_B\rangle)$;
\item [(3)] $\lim_{n\to\infty}\Vert \phi_n(x)-x\Vert_2=0$ for every $x\in N$.
\end{enumerate}
\end{Def}
As for rigidity, relation property H is independent of the chosen trace.

\bigskip

Typical and important examples of inclusions of finite von Neumann algebras with property H as above come from crossed products: assume that a countable group $G$ acts on some finite von Neumann algebra $B$ and that the action, denoted by $\sigma$, preserves some trace $\tau$.  

Then one has:

%Prop 4.2
\begin{Prop}\label{ExRelH}
(\cite{PopaAnnals}, Proposition 3.1) Let $G$, $(B,\tau)$ and $\sigma$ be as above. Then $N=B\rtimes_\sigma G$ has property H relative to $B$ if and only if the group $G$ has the Haagerup property. \eprsk
\end{Prop}

Now, let us consider the following situation: let $\Gamma$ be a countable group that acts by automorphisms on some discrete group $H$. Then the von Neumann algebra $L(H\rtimes\Gamma)$ associated to the semidirect product group $G=H\rtimes \Gamma$ is naturally isomorphic to the crossed product algebra $N=L(H)\rtimes \Gamma$. Assume that $\Gamma$ has the Haagerup property and that the pair $(G,H)$ has the relative property (T). Then $L(H)$ is an HT subalgebra of $N$ in the sense of the following definition:

%Def 4.3
\begin{Def}\label{DefHT}
(Definition 6.1, \cite{PopaAnnals})
Let $N$ be a finite von Neumann algebra and let $B\subset N$ be a von Neumann subalgebra. Then $B$ is an \textit{HT subalgebra} of $N$ if the following two conditions are met:
\begin{enumerate}
\item [(a)] $N$ has property H relative to $B$;
\item [(b)] there exists a von Neumann subalgebra $B_0\subset B$ such that $B'_0\cap N\subset B$ and the inclusion $B_0\subset N$ is rigid. 
\end{enumerate}
If $A\subset N$ is a Cartan subalgebra that satisfies both conditions, then it is called an \textit{HT Cartan subalgebra}.
\end{Def}

The importance of existence of HT Cartan subalgebras is illustrated by the following results:

%Thm 4.4
\begin{Thm}\label{ThmHTCartan}
(Theorem 6.2, \cite{PopaAnnals}) Let $N$ be a type $\mathrm{II}_1$ factor with two HT Cartan subalgebras $A_1$ and $A_2$. Then there exists a unitary element $u\in N$ such that $uA_1u^*=A_2$. \eprsk
\end{Thm}

%Cor 4.5
\begin{Cor}\label{CorHTCartan}
(Corollaries 6.5 and 8.2, \cite{PopaAnnals})
Let $A_i\subset N_i$, $i=1,2$, be HT Cartan subalgebras of type $\mathrm{II}_1$ factors, and let $\theta:N_1\ra N_2$ be an isomorphism. Then there exists $u\in U(N_2)$ such that $u\theta(A_1)u^*=A_2$. Thus, if $A\subset N$ is an HT Cartan subalgebra of a type $\mathrm{II}_1$ factor $N$,
there exists a unique (up to isomorphism) standard equivalence relation $\Rr_N^{HT}$ on the standard probability space $(X,\mu)$ such that the inclusion $A\subset N$ is isomorphic to the inclusion $L^\infty(X,\mu)\subset L(\Rr_N^{HT},s)$ for some cocycle $s$. Thus the factor $N$ inherits $\ell^2$-Betti numbers $\beta_k^{HT}(N)\coloneqq \beta_k(\Rr_N^{HT})$ for every integer $k\geq 0$, which are isomorphism invariants of $N$. \eprsk
\end{Cor}

Following Notation 6.4 of \cite{PopaAnnals}, we denote by $\mathcal{HT}$ the class of all type II$_1$ factors with separable predual with HT Cartan subalgebras.

Another interesting property of the $\ell^2$-Betti numbers $(\beta_k^{HT}(N))_{k\geq 0}$ for $N\in \mathcal{HT}$ is that, for every real number $t>0$, we have $\beta_k^{HT}(N^t)=\beta_k^{HT}(N)/t$ (Corollary 8.2, \cite{PopaAnnals}), where $N^t$ denotes the \textit{amplification of $N$ by $t$}, which is the factor $pM_n(N)p$, where $n\geq t$ is a positive integer and $p\in M_n(N)$ is a projection with (non-normalized) trace $\Tr\otimes\tau(p)=t$. (We recall that the isomorphism class of $N^t$ depends neither on $p$ nor on $n\geq t$.)

Recall also that if $N$ is a type II$_1$ factor, its \textit{fundamental group} $\mathcal{F}(N)$ is the set of all $t>0$ such that $N^t$ is isomorphic to $N$. Hence, if $N\in\mathcal{HT}$ is such that $0<\beta_k(N)<+\infty$ for some $k$, then $\mathcal{F}(N)$ is trivial.

\par\vspace{2mm}

{\emph{Proof of Theorem \ref{ThmTrivialFundGr}:}}
Set $A_n=L^\infty(\T^n)$ which is a Cartan subalgebra of $N_n$. Then, by Proposition 3.1 of \cite{ValMax},  the inclusion $A_n\subset N_n$ is rigid, and as $G_n$ has the Haagerup property, $A_n$ is an HT Cartan subalgebra of $N_n$. Furthermore, it follows from Corollary 3.16 in \cite{Ga}, that $\beta_k(\Rr_N^{HT})=\beta_k(G_n)$ for all $k\geq 0$ and $n\geq 2$. The rest is a consequence of the above discussion and Corollary \ref{CorHTCartan}.
\eprsk

%Remark 4.6
\begin{Rem} Theorem \ref{ThmTrivialFundGr} raises a natural question: are $N_m$ and $N_n$ isomorphic or not for $m\neq n$? It was pointed out to us by A. Ioana that, by Theorem 1 in his paper \cite{Ioa12}, the factor $N_n$ has a unique group measure space Cartan subalgebra, up to unitary conjugation. As a consequence, $N_m$ is isomorphic to $N_n$ if and only if the equivalence relation $\mathcal{S}_m$ on $\T^m$ is isomorphic to the equivalence relation $\mathcal{S}_n$ on $\T^n$, or, as both are orbital with respect to free and ergodic actions, if and only if the latter are orbit equivalent.
\end{Rem}

%Ssection 4.2
\subsection{PROOF OF THEOREM \ref{MaxHaagerup}}

First, we need to prove the existence of maximal Haagerup von Neumann algebras in finite ones.

In fact, we will prove a more general result which not only yields the above statement, but which implies the following result of independent interest: given a finite von Neumann algebra $(N,\tau)$ with separable predual, the subset $\HAP(N)$ of all von Neumann subalgebras of $N$ that have the Haagerup property is closed with respect to the \textit{Effros-Mar\'echal topology} that we recall now (see for instance \cite{HW} or \cite{AHO}). Denote by $\SA(N)$ the set of all von Neumann subalgebras of $N$. As the latter is finite, the topology on $\SA(N)$ can be described as follows: a net $(M_i)_{i\in I}\subset \SA(N)$ converges to $M\subset N$ if and only  
\[
\liminf_{i\in I} M_i=\limsup_{i\in I} M_i=M,
\] 
where
\[
\liminf_{i\in I} M_i\coloneqq\{x\in N\colon \exists (x_i)_{i\in I}\in\ell^\infty (I,M_i)\ :\  \lim_{i\in I}\Vert x_i-x\Vert_2=0\}
\]
and 
\[
\limsup_{i\in I} M_i\coloneqq \la \{x\in N\colon \exists (x_i)_{i\in I}\in \ell^2(I,M_i)\ : \ \mathrm{wo}-\lim_{i\in I} x_i=x\}\rg.
\] 
Notation $\la E\rg$ means the von Neumann algebra generated by the subset $E\subset N$, wo-$\lim$ refers to the limit with respect to the weak operator topology on $N$ and $\ell^\infty (I,M_i)$ denotes the C$^*$-algebra of bounded nets $(x_i)_{i\in I}\in\prod_i M_i$.

\bigskip
Statement (1) of the next result was inspired by Proposition 4.7 in \cite{AHO} that states that the subset of amenable von Neumann subalgebras of $N$ is closed; it is a von Neumann algebraic analogue of the fact that, for a countable group $G$, the set of Haagerup subgroups is closed in the Chabauty topology on the set of all subgroups of $G$ (which follows easily from Proposition 6.1.1 in \cite{CCJJV}). 

%Prop 4.7
\begin{Prop}\label{ExistHMax}
Let $(N,\tau)$ be a finite von Neumann algebra as above and let $(M_i)_{i\in I}\subset \HAP(N)$ be a net such that $M\coloneqq \liminf_{i\in I}M_i$ exists. Then $M\in \HAP(N)$. In particular, 
\begin{enumerate}
\item [(1)] the subset $\HAP(N)$ is closed in $\SA(N)$ with respect to the Effros-Mar\'echal topology and,
\item [(2)] for every von Neumann subalgebra $P\in\HAP(N)$, there exists a maximal Haagerup von Neumann subalgebra $Q\subset N$ which contains $P$.
\end{enumerate}
\end{Prop}
\bpr Denoting by $(N)_1$ the unit ball of $N$ for the operator norm, and setting $(E)_1=E\cap (N)_1$ for every subset $E$ of $N$, let $(M_i)_{i\in I}\in \HAP(N)$ and $M\in \SA(N)$ be such that $M=\liminf_{i\in I}M_i$.\\
To show that $M\in \HAP(N)$, it suffices to prove that, for every finite set $F\subset (M)_1$ and every $\eps>0$, there exists a completely positive map $\Phi=\Phi_{(F,\eps)}:M\ra M$ with the following properties:
\begin{enumerate}
\item [(a)] $\Phi(1)\leq 1$ and $\tau\circ \Phi\leq \tau$;
\item [(b)] the extension $T_\Phi$ of $\Phi$ on $L^2(M,\tau)$ is a compact operator;
\item [(c)] $\max_{x\in F}\Vert \Phi(x)-x\Vert_2\leq \eps$.
\end{enumerate}
Thus, let $F\subset (M)_1$ and $\eps>0$ as above. Let us number the elements of $F$ so that $F=\{x_1,\ldots,x_r\}$. For each $1\leq k\leq r$, there exists a net $(x_k^{(i)})_{i\in I}\in \ell^\infty(I,M_i)$ such that $x_k^{(i)}\in (M_i)_1$ for every $i$ and that
\[
\lim_{i\in I}\Vert x_k^{(i)}-x_k\Vert_2=0.
\] 
Fix henceforth an element $j\in I$ such that
\[
\max_{1\leq k\leq r}\Vert x_k^{(j)}-x_k\Vert_2\leq \eps/3.
\]
Since $M_j\in \HAP(N)$, there exists a completely positive map $\phi_j:M_j\ra M_j$ which has properties (a) and (b), and such that 
\[
\max_{1\leq k\leq r}\Vert\phi_j(x_k^{(j)})-x_k^{(j)}\Vert_2\leq \eps/3.
\]
Set $\Phi\coloneqq E_M\circ \phi_j\circ E_{M_j}|_M: M\ra M$.
It obviously satisfies condition (a) above. Moreover, $T_\Phi=e_MT_{\phi_j}e_{M_j}|_{L^2(M,\tau)}$ is a compact operator on $L^2(M,\tau)$, and finally we have for every $1\leq k\leq r$, since $x_k\in M$,
\begin{align*}
\Vert \Phi(x_k)-x_k\Vert_2
&=
\Vert E_M[\phi_j(E_{M_j}(x_k))-x_k]\Vert_2
\leq 
\Vert \phi_j(E_{M_j}(x_k))-x_k\Vert_2\\
&\leq
\Vert \phi_j[E_{M_j}(x_k)-x_k^{(j)}]\Vert_2+\Vert\phi_j(x_k^{(j)})-x_k^{(j)}\Vert_2+\Vert x_k^{(j)}-x_k\Vert_2\\
&\leq
2\Vert x_k^{(j)}-x_k\Vert_2+\Vert \phi_i(x_k^{(j)})-x_k^{(j)}\Vert_2\leq \eps.
\end{align*}
(We have used the inequality $\Vert\phi_j(y)\Vert_2\leq \Vert y\Vert_2$ for every $y\in M_j$ which follows from $\tau\circ\phi_j\leq \tau$.)\\
Statement (1) follows then trivially, and, in order to prove Statement (2), 
by Zorn's lemma, it suffices to prove that, if $(Q_i)_{i\in I}$ is an increasing net in $\HAP(N)$ such that $P\subset Q_i$ for every $i$, then $Q\coloneqq (\bigcup_{i\in I}Q_i)''\in \HAP(N)$. But obviously $Q=\liminf_i Q_i$, and $Q\in \HAP(N)$ by the main part of the proof.
\eprsk

\bigskip
{\emph{Proof of Theorem \ref{MaxHaagerup}:}}
Fix $n\geq 2$. The abelian subalgebra $A_n=L^\infty(\T^n)$ being a Cartan subalgebra of $N_n$, by Proposition 3.6 of \cite{PopaInv}, it is contained in an intermediate hyperfinite subfactor $P$ of $N_n$ in which it is still a Cartan subalgebra. Let then $P\subset Q\subset N_n$ be a maximal Haagerup von Neumann subalgebra whose existence is guaranteed by Proposition \ref{ExistHMax}. By the considerations following Definition \ref{DefDicho}, we see that $Q$ is a factor, hence $Q=L(\mathcal{R},s)$, where $\mathcal{R}$ is an ergodic subequivalence relation of $\mathcal{S}_n$ and $s$ is a Borel 2-cocycle on $\mathbb{R}$. By Corollary \ref{DichoWithCoc}, either $\mathcal{R}$ is amenable or the inclusion $L^\infty(X,\mu)\subset L(\Rr,s)=Q$ is rigid. The second case is precluded by the assumption that $Q$ has the Haagerup property, so $\mathcal{R}$ is amenable, hence $s$ is trivial and $Q$ is a hyperfinite von Neumann subalgebra. 

Next, fix $n\geq 2$ even. Let $H$ be a maximal amenable subgroup of $\SL_2(\Z)$ containing some hyperbolic matrix. Set $P=A_n\rtimes H$; view $P$ as a von Neumann subalgebra of $N_n$. Then $P$ is a subfactor, by the second part of Theorem \ref{ergodic+free}. Moreover $P$ is maximal amenable in $N_n$, by a result of Boutonnet and Carderi (Theorem A and Corollary B in \cite{BoCa}). The above argument shows that $P$ is maximal Haagerup in $N_n$.

Finally, if $n\geq 3$ is odd, and if $P$ and $Q$ are as in the first part of the proof, being factors, they cannot be of the form $L(K)$ for any subgroup $K$ of $\Z^n\rtimes \PSL_2(\Z)$ by Proposition \ref{nonergodic}.
\eprsk

\bigskip
\noindent
{\bf Authors'address}:\\
Institut de Math\'ematiques, Universit\'e de Neuch\^atel, Rue E.-Argand 11, 2000 Neuch\^atel, Switzerland.\\
pajolissaint@gmail.com, alain.valette@unine.ch

\bigskip
\noindent
{\bf Acknowledgements:} We thank C. Anantharaman-Delaroche, J. Blanc, D. Gaboriau, A. Ioana, M. Kalantar, N. Monod, J. Renault and S. Vaes for useful exchanges. The referee's careful reading allowed to eliminate a gap in the original proof of Theorem \ref{IoanaDicho}.
%\end{acknowledgements}

%\bibliographystyle{plain}
%\bibliography{refFWap}

\end{document}